\documentclass{cimart}

\usepackage[capitalise]{cleveref}

\newtheorem{prob}{Problem}

\newcommand{\alphaa}[2]{\alpha_{#1}{#2}}
\newcommand{\betaa}[2]{\beta_{#1}{#2}}
\newcommand{\thetaa}[2]{\theta_{#1}{#2}}

\DeclareMathOperator{\A}{A}
\DeclareMathOperator{\E}{E}
\DeclareMathOperator{\End}{End}
\DeclareMathOperator{\Ext}{Ext}
\DeclareMathOperator{\id}{id}
\DeclareMathOperator{\im}{Im}
\DeclareMathOperator{\Ker}{Ker}
\DeclareMathOperator{\Ret}{Ret}
\DeclareMathOperator{\sgn}{sgn}
\DeclareMathOperator{\Sym}{Sym}
\DeclareMathOperator{\tr}{tr}

\VOLUME{33}
\ISSUE{3}
\YEAR{2025}
\NUMBER{2}
\DOI{https://doi.org/10.46298/cm.16958}

\title{Set-theoretical solutions to the pentagon equation: a survey}

\author{Marzia Mazzotta}

\authorinfo{Universit\`a del Salento, Italy}{marzia.mazzotta@unisalento.it}

\abstract{This survey aims to collect the main results of the theory of the set-theoretical solutions to the pentagon equation obtained up to now in the literature. In particular, we present some classes of solutions and raise some questions.}



\begin{document}

\section*{Introduction}
The \emph{pentagon equation} classically originates from the field of Mathematical Physics. It appears in various forms in different contexts, such as in the representation theory of quantum groups as the Biedenharn–Elliott identity for $6j$-symbols \cite{BiLo81}, in quantum conformal field theory as the identity for the fusion matrices \cite{MoSe89}, in quasi-Hopf algebras as the consistency equation for the association \cite{Dr89}. Still, it appears in other areas of mathematics, also with different terminologies (see, for instance, \cite{ BaSk93, BaSk98, Ka23, Ka11, Kawa10, Mi98, Wo96}). The pentagon equation belongs to the family of \emph{polygon equations} \cite{DiMu15} that are associated with a Tamari lattice, namely, a partially ordered set whose elements consist of all the binary bracketings of a string of $n+1$ symbols. In particular, the pentagon equation describes a Tamari lattice of order four.  We refer to the papers \cite{DiMu15, Mu23x} for more details on this topic and further references to the occurrences of the pentagon equation in literature. 

Given a vector space $V$ over a field $F$, a linear map $S: V \otimes V \to V \otimes V$ is said to be a \emph{solution of the pentagon equation}  if 
	\begin{align*}
	S_{12}S_{13}S_{23}=S_{23}S_{12},
	\end{align*}
	where 
		$S_{12}=S\otimes \id_V$,\,$S_{23}=\id_V\otimes \, S$,\,$S_{13}=(\id_V\otimes\, \Sigma)\,S_{12}\;(\id_V\otimes \, \Sigma)$,
	with $\Sigma$ the flip operator on $V\otimes V$, i.e., $\Sigma(u\otimes v)=v\otimes u$, for all $u,v\in V$. Maillet \cite{Ma94} showed that solutions of the pentagon equation lead to solutions of the tetrahedron equation \cite{Za80}, a multidimensional generalization of the well-known quantum Yang-Baxter equation \cite{Ba72, Ya67}. Indeed, as one can note, the pentagon equation is the quantum Yang-Baxter equation with the middle term missing on the right-hand side. By the way, we highlight that the Yang-Baxter equation and the tetrahedron equation belong to the family of \emph{$N$-simplex equations} which are related to the higher Bruhat order (see \cite{DiMu15}).  Furthermore, Militaru showed that bijective solutions on finite dimensional vector spaces are in one-to-one correspondence with finite dimensional Hopf algebras, so the classification of the latter is reduced to the classification of solutions (see \cite[Theorem 3.2]{Mi04}). 
 
In 1998, Kashaev and Sergeev \cite{KaSe98} began the study of the pentagon equation in set-theoretical terms to find vector solutions, pursued further in
\cite{KaRe07}. If $X$ is a finite set and $s:X\times X\to X\times X$ is a map satisfying the relation
 \begin{equation*}
	s_{23}s_{13}s_{12}=s_{12}s_{23},
	\end{equation*}
	where
		$s_{12}=s\times \id_X$,  $s_{23}=\id_X\times s$, $s_{13}=(\id_X\times \tau)\,s_{12}\,(\id_X\times \tau)$, with $\tau(x,y)=(y,x)$, for all $x,y\in X$,  they associated  its pullback $S$, i.e., the linear operator in the space of the functions on $X \times X$ defined as $S(f)(x,y)= f(s(x,y)),$ for all $x, y \in X$.
The map $s$ above is said to be a \emph{set-theoretical solution of the pentagon equation}, or briefly a \emph{PE solution}, on $X$. Thus, if $s$ is a PE solution on $X$, then $S$ is a solution of the pentagon equation with appropriate definition for the tensor product of infinite dimensional vector spaces. Yet, the first instances of PE solutions can be also extrapolated from the papers by Zakrzewski \cite{Za92} 
and Baaj and Skandalis \cite{BaSk98}, where the authors dealt with specific solutions on differential manifolds and measure spaces, respectively. A purely algebraic transcription of these solutions can be found in \cite{CaMaMi19}. In this paper, the authors provide the first systematic approach for investigating PE solutions. Writing a PE solution $s: X \times X \to X \times X$ as $s(x,y)=(x \cdot y, \theta_x(y))$, where $\theta_x:X\to X$ is a map, for every $x \in X$,  one has that 
\begin{itemize}
    \item[-] $(X, \cdot)$ is a  semigroup,
    \item[-] $\theta_x(y) \cdot  \theta_{x \cdot y}(z)=\theta_x(y \cdot z )$,
    \item[-] $\theta_{\theta_x(y)}\theta_{x\cdot y}=\theta_y$,
    \end{itemize}
    for all $x,y,z \in X$. 
  Therefore, the natural first step is to look for PE solutions in a group. In this regard, a complete classification of these maps was provided in \cite[Theorem 13]{CaMaMi19} in terms of normal subgroups of the group. Subsequently, the authors in \cite{MaPeSt23x} have started the study of PE solutions in the family of Clifford semigroups with an attempt to describe them (see \cite{Cl41} and the monographs \cite{Pe84, La98} for more details on this class of semigroups). However, the study has not been completed because there are many families of solutions.  Given a Clifford semigroup $X$, they describe the \emph{$\E(X)$-invariant PE solutions}, namely, those solutions for which $\theta_a(e)=\theta_a(f)$, for all $e,f \in\E(X)$ and for every $a\in X$. 
Moreover, they construct a family of \emph{$\E(X)$-fixed PE solutions},  i.e., those solutions for which $\theta_a$ fixes every element in $\E(X)$, for every $a\in X$. Further developments in this regard can be found in this paper.

Describing PE solutions on arbitrary semigroups seems very challenging since there are many even in the case of small-order semigroups. For example, for semigroups of order $3$ there are already 202 non-isomorphic solutions (see \cite[Appendix B]{tesi}).
Because of this, several authors have studied specific classes of PE solutions.  In this respect, Colazzo, Jespers, and Kubat \cite{CoJeKu20} described all the \emph{involutive} PE solutions, 
namely, solutions satisfying the property $s^2=\id_{X \times X}$. Recently, the author in \cite{Maz23} started the investigation of \emph{idempotent} PE solutions, i.e., solutions such that $s^2=s$. A classification of these solutions is obtained in the case of monoids having central idempotents. Additionally, in \cite{CaMaSt20}, one can find a description of PE solutions that are also set-theoretical solutions of the Yang-Baxter equation, the so-called \emph{P-QYBE solutions}. Easy examples are the maps $s$ of the form $s(x, y) = (f(x), g(y))$, with $f, g: X \to X$ idempotent maps such that $fg = gf$, that belong to the class of Lyubashenko's solutions (for more details, see \cite{Dr92}).

Another line of research consists of looking for methods to obtain new solutions starting from known ones. In \cite{CaMaSt20}, some can be found, such as the technique of the \emph{matched product} of two solutions $s$ and $t$ on two
semigroups $S$ and $T$, respectively, that is a new solution on the semigroup $S \bowtie T$ (see \cite[Definition 1]{CaMaSt20}), namely a semigroup including the \emph{classical Zappa product} of $S$ and $T$ (see \cite[Definition 1.1]{Ku83}).

In this survey, we give an overview of PE solutions and collect several results on this little-explored area. Some open questions will be raised during the work that may inspire the reader. Specifically, the survey is organized as follows: in the first section, we give some preliminaries and examples by recalling the description of all PE solutions defined in groups; in the second one, we present some problems related to PE solutions in the class of Clifford semigroups; the third section is devoted to the description of all involutive PE solutions: the fourth one contains the characterization of idempotent PE solutions in monoids having central idempotents and some issues concerning this class of solutions. Finally, in the last paragraph, we give some hints on other classes of PE solutions that one could study, namely the commutative and cocommutative ones.

\section{Definitions and examples}
This section aims to give some basics on PE solutions. In particular, we recall some examples and constructions.

Let $X$ be a non-empty set. The map $s: X \times X \to X \times X$ given by $s(x,y)=(x \cdot y, \theta_x(y))$, where $\theta_x:X\to X$ is a map, for every $x \in X$,  is a \emph{PE solution} if, and only if,  $(X, \cdot)$ is a semigroup and the following hold
	\begin{align}
	\label{p_one}\theta_x(y) \cdot  \theta_{x \cdot y}(z)&=\theta_x(y \cdot z ), \tag{P1}\\
	\label{p_two}\theta_{\theta_x(y)}\theta_{x\cdot y}&=\theta_y, \tag{P2}
	\end{align}
	for all $x,y,z \in X$. Moreover, $s$ is said to be \emph{involutive} if $s^2=\id_{X \times X}$; \emph{idempotent} if $s^2=s$; \emph{non-degenerate} if $\theta_x$ is bijective, for every $x \in X$. 
 \begin{definition}
    Let $X$ and $Y$ be two semigroups and consider $s(x,y)=(xy, \theta_x(y))$ and $t(u, v)=(uv, \eta_u(v))$ two PE solutions on $X$ and $Y$, respectively. Then, $s$ and $t$ are \emph{isomorphic} if there exists a semigroup isomorphism $\psi: X \to Y$ such that $\psi\theta_x(y)=\nolinebreak\eta_{\psi\left(x\right)}\psi(y)$,
    for all $x,y \in X$, or, equivalently, $(\psi \times \psi)s=t( \psi\times \psi)$.    
    \end{definition}
 From now on, we will use
juxtaposition for the product in any semigroup.
In the following, we give some examples of solutions.
\begin{example} \label{exsKey}\hspace{1mm}
    \begin{enumerate}
        \item \label{ex_militaru} Let $X$ be a set and $f, g: X \to X$ maps such that $f^2=f$, $g^2=g$, and $fg=gf$. Then, the map 
	$s(x,y)=\left(f\left(x\right),\, g\left(y\right)\right)$, for all $x, y \in X$, is an idempotent PE solution on $X$ (cf. \cite{Mi98}). Note that this map also belongs to the class of P-QYBE solutions since it also satisfies the quantum Yang-Baxter equation (cf. \cite[Examples 8]{CaMaSt20}).
 \item If $S$ is a semigroup and $\gamma\in \End(S)$ such that $\gamma^2=\gamma$, the map 
	$s(x,y)=\left(xy,\gamma\left(y\right)\right),$
	for all $x,y \in S$, is a PE solution on $S$ (see \cite[Examples 2-2.]{CaMaMi19}). Note that $s$ is non-degenerate if, and only if, $\gamma=\id_S$.
 \item  Let $S=\{0,a,b\}$  be the null semigroup, i.e., $xy=0$, for all $x, y \in S$. Consider the maps $\theta_0=\id_S$ and $\theta_a=\theta_b$ such that $\theta_a(0)=0$, $\theta_a(a)=b$, and $\theta_a(b)=a$. Then, the map $s(x, y)=(0, \theta_x(y))$ is an idempotent and non-degenerate PE solution on $S$ (see \cite[Appendix B]{tesi}).
 \item Let $G$ be a group of finite exponent, $E=\{1, \dots, n\}$, and $\sigma \in \Sym(n)$ such that $\sigma^{\sigma(i)+1}=\sigma^i$, for any $i \in E$. Put $S=G \times E$, the map $s: S \times S \to S \times S$ defined by
 \begin{align*}
     s((i,a), (j, b))=\left(\left(i,ab\right), \left(\sigma^i(j), b\right)\right)
 \end{align*}
 is a bijective PE solution on $S$ (see \cite[Example 1.3]{CoJeKu20}). Additionally, $s$ is involutive if, and only if $G$ has exponent $2$ and $\sigma$ has order $2$.
    \end{enumerate}
\end{example}

Other examples of PE solutions can be obtained from the pioneering works \cite{Za92, BaSk98}, where the authors deal with the PE solutions in non-algebraic contexts, by calling such maps \emph{pentagonal transformations}. We present such solutions by transcribing them in purely algebraic terms.

\begin{example} \label{exs_za_bask}	Let $G$ be a group that admits an exact factorization through two subgroups $H,K$ and 
	$p_1:G \to H$, $p_2: G \to K$ the projection maps such that every $x\in G$ can be written as $x=p_1(x)\,p_2(x)$. Then, the maps $s, r:G\times G\to G\times G$ given by 
\begin{align*}
s\left(x,y\right)=\left(p_2\left(yp_1\left(x\right)^{-1}\right)x,\; yp_1\left(x\right)^{-1}\right) \quad \text{and} \quad 
r\left(x,y\right)=\left(xp_1\left(p_2\left(x\right)^{-1}y\right), \, p_2\left(x\right)^{-1}y\right)   
\end{align*}
are bijective PE solutions on $G$. In particular, $G$ endowed with the operation defined by the first component of the map $r$ (or $s$) is a left group, i.e., the direct product of a left zero semigroup and a group.
\end{example}

The maps in \cref{exs_za_bask} are such that $r= \tau s \tau$, with $\tau$ the flip map.
In general, one has that if $s$ is a bijective PE solution on a set $X$, the map $s^{op}=\tau s^{-1}\tau$ also is a PE solution, called the \emph{opposite PE solution} of $s$ (in analogy to the opposite operators defined on Hilbert spaces \cite{BaSk98}).

In \cite[Proposition 1]{KaSe98}, one can find the first systematic way for constructing PE solutions which we recall below.

\begin{proposition}\label{prop_kase_costr}
		Let $G$ be a group, $X$ a subsemigroup of $G$, and $\lambda, \mu : X \to G$ maps such that, for all $x, y \in X$,
  \begin{enumerate}
      \item $x \ast y:=\mu(x)^{-1} \mu(x y) \in X$,
      \item $\mu(x \ast y)=\lambda(x)\mu(y)$.
  \end{enumerate}
Then, the map $s: X \times X \to X \times X$ given by $s(x,y)=(x y, x \ast y)$ is a PE solution on $X$.
  \end{proposition}

The following theorem illustrates how to find all the PE solutions in the case of a group. 
	\begin{theorem}\cite[Theorem 15]{CaMaMi19}\label{teo_gruppi}
	Let $G$ be a group, $K$ a normal subgroup of $G$, and $R$ a system of representatives of $G/K$.  If $\mu: G \to R$ is a map such that $\mu(x) \in Kx$, for any $x \in G$, then the map $
	s(x,y)=\left(xy, \mu\left(x\right)^{-1}\mu\left(xy\right)\right),
$
	for all $x,y \in G$, is a PE solution. 
 
	Vice versa, if $s(x,y)=(xy, \theta_x(y))$ is a PE solution on the group $G$, then the set
	$K=\{x \in G\, \mid\, \theta_1(x)=1\}$
	is a normal subgroup of $G$, $\im \theta_1$ is a system of representatives of $G/K$, $\theta_1(x) \in Kx$, and $\theta_x(y)=
	 \theta_1\left(x\right)^{-1}\theta_1\left(xy\right),$
	for all $x, y \in G$.
\end{theorem}

  It immediately follows from \cref{teo_gruppi} that the only bijective PE solution on a group $G$ is given by the map $s_G(x, y)=(xy, y)$ that is known in literature as \emph{Kac-Takesaki solution}, concerning the associated linear operator defined on the Hilbert space of square-integrable functions \cite{Tak72}.

An easy example of PE solutions on the symmetric group is the following.

\begin{example} Let $n \geq 3$, $\Sym_n$ the symmetric group of order $n$, $\A_n$ the alternating group of degree $n$, and $R=\lbrace \id_{\Sym_n}, \pi \rbrace$ a system of representatives of $\Sym_n/\A_n$, where $\pi$ is a transposition in $\Sym_n$. Moreover, consider the map $\mu: \Sym_n \to R$ given by
	\begin{center}
		$\mu(\alpha)=\begin{cases}
		\pi \qquad &\text{if} \;\sgn(\alpha)=-1 \\
		\id_{\Sym_n} \quad&\text{if}\; \sgn(\alpha)=1,
		\end{cases}$
	\end{center}
	for every $\alpha \in \Sym_n$. Then, the map $s(\alpha,\beta)=\left(\alpha\beta,\mu\left(\alpha\right)^{-1}\mu\left(\alpha\beta\right)\right)$, for all $\alpha, \beta \in \Sym_n$, is a PE solution on $\Sym_n$.
\end{example}

In the next section, we will show how to extend the theorem above to a specific class of solutions defined on Clifford semigroups. Regarding the PE solutions defined on monoids, some properties of the maps $\theta_x$ are given in \cite{Maz23} and we recall them below. 

\begin{lemma}\label{prop_monoid}
Let $M$ be a monoid with identity $1$ and $s(x,y)=(xy, \theta_x(y))$ a PE solution on $M$. Then, the following hold:
\begin{enumerate}
    \item $\theta_x(1) \in \E(M)$,
    \item $\theta_1=\theta_{\theta_x(1)}\theta_x$,
    \item $\theta_1(x) \in \theta_1(1)M$,
    \item $\theta_x=\theta_{\theta_1(x)}\theta_x$,   
    \end{enumerate}
    for every $x \in M$.
\end{lemma}

However, we are still far from a classification.

\begin{prob}
Study or find some constructions of PE solutions defined on monoids.
\end{prob}

To conclude this preliminary section, we mention the technique of \emph{matched product} of PE solutions introduced in \cite{CaMaSt20}, which is a method to obtain a new PE solution on the Cartesian product of two semigroups starting from two given ones. To recall it, we first give the following definition.
\begin{definition}
	Let $S$ and $T$ be semigroups, $\alpha:T\to S^{S}$ and $\beta:S\to T^{T}$ maps, and set $\alpha_{u}:= \alpha\left(u\right)$, for every $u\in T$, and $\beta_{a}:= \beta\left(a\right)$, for every $a\in S$. Moreover, let $s(a,b)=(ab, \theta_a(b))$ and $t(u,v)=(uv, \eta_u(v))$ PE solutions on $S$ and $T$, respectively. If the following conditions are satisfied
	\begin{align*}
	\alpha_{u}\left(a\alpha_{v}\left(b\right)\right)
	= \alpha_{u}\left(a\right)\alpha_{\beta_{a}\left(u\right)v}\left(b\right) \quad &\text{and} \quad
	\beta_{a}\left(\beta_{b}\left(u\right)v\right)
	= \beta_{b\alpha_{v}\left(a\right)}\left(u\right)\beta_{a}\left(v\right), \\
	\thetaa{a}{\alphaa{u}{}} &= \thetaa{\alphaa{v}{\left(a\right)}}{\alphaa{\betaa{a}{\left(v\right) u}}{}}, \\ \thetaa{a\alphaa{u}{\left(b\right)}}{}&=\alphaa{\eta_{\betaa{b}{\left(u\right)}}{\left(v\right)}}{\thetaa{a\alphaa{u}{\left(b\right)}}{}}, \\
	\eta_{\betaa{b\alphaa{v}{\left(c\right)}}{\left(u\right)}}{\betaa{c}{\left(v\right)}}&=\betaa{\thetaa{a\alphaa{u}{\left(b\right)}}{\alphaa{\betaa{b}{\left(u\right)v}}{\left(c\right)}}}{\eta_{\betaa{b}{\left(u\right)}}{\left(v\right)}} ,
	\end{align*}
	for all $a,b,c\in S$ and $u,v\in T$, then we call $\left(s,t,\alpha,\beta\right)$ a \emph{matched quadruple}.
\end{definition}

  It is a routine computation to verify that $S \times T$, endowed with the operation defined by $
	\left(a,u\right)\left(b,v\right)=\left(a\alphaa{u}{\left(b\right)}, \betaa{b}{\left(u\right)}v\right)
$, for all $(a,u), (b,v) \in S \times T$,	is a semigroup, called  \emph{matched product of $S$ and $T$} and denoted by $S \bowtie T$. As a class of examples, one can easily find the classical Zappa product \cite{Ku83}. 

\begin{theorem}\label{TeoMatch}
	Let $S$, $T$ be semigroups and $\left(s,t,\alpha,\beta\right)$ a matched quadruple. Then, the map $s \bowtie t:(S \times T) \times (S \times T) \to (S \times T) \times (S \times T)$ defined by
	\begin{align*}
	s \bowtie t\left((a,u) , (b,v)\right)
	= \left(\left(a\alphaa{u}{\left(b\right)}, \betaa{b}{\left(u\right)}v\right), \,  \left(\thetaa{a}{\alphaa{u}{\left(b\right)}}, \eta_{\betaa{b}{\left(u\right)}}{\left(v\right)}\right)\right), 
	\end{align*}
	for all	$\left(a,u\right),\left(b,v\right)\in S\times T$, is a PE solution on $S \bowtie T$, called the \emph{matched product of $s$ and $t$}.
 \end{theorem}

\begin{example} \label{ex-match-quadr}
Let $S=\lbrace 1_S,x,y\rbrace$ be the commutative idempotent monoid such that $xy=y$ and $T=\lbrace 1_T, z\rbrace$ the unique idempotent monoid on two elements. Consider the maps $\alpha: T \to S^S$ given by $\alpha_{1_T}=\id_S$ and $\alpha_z=\gamma$, where $\gamma \in \End(S)$, $\gamma^2=\gamma$, and $\beta: S \to T^T$ defined by $\beta_{1_S}=\id_T$ and $\beta_x(u)=\beta_y(u)=1_T$, for every $u \in T$.
Moreover, let $s(a,b)=(ab, \gamma(b))$ and $t(u,v)=(uv,v)$ be the PE solutions on $S$ and $T$, respectively. Then, $(s,t,\alpha, \beta)$ is a matched quadruple and the map
    \begin{center}
		$s \bowtie t \, ((a,u) , (b,v))=\begin{cases}
		((a,v), ( 1_S,v)) \qquad &\text{if} \;b=1_S, \,u=1_T \\
		((a,uv) , (1_S,v)) \qquad &\text{if} \;b= 1_S, \,u\neq 1_T\\
		((ab,v, (\gamma(b),v) )\qquad &\text{if} \;b\neq 1_S, \,u=1_T\\
		((a\gamma(b),v), (\gamma(b),v) )\qquad &\text{if} \;b\neq 1_S, \,u\neq 1_T
		\end{cases}$ 
	\end{center}
is a PE solution on $S \bowtie T$.
\end{example}

\section{PE solutions on Clifford semigroups}
In this section,  we focus on PE solutions defined on a Clifford semigroup and satisfying special properties on the set of all idempotents discussed in \cite{MaPeSt23x}. Moreover, we propose some topics to be developed.

We begin by recalling some basic facts
from the theory of Clifford semigroups that will be useful in
the following. For a more detailed treatment of this class of semigroups, we refer the reader to the paper by Clifford \cite{Cl41} and the monographs by Petrich \cite{Pe84} and Lawson \cite{La98}.
A semigroup $X$ is \emph{inverse} if for each $x\in X$ there exists a unique $x^{-1}\in X$ such that \begin{center}
    $x = xx^{-1}x$ \qquad \text{and}\qquad $x^{-1} = x^{-1}xx^{-1}$.
\end{center} 
Clearly, every group is an inverse semigroup. We have $(xy)^{-1}=y^{-1} x^{-1}$ and $(x^{-1})^{-1}=x$, for all $x,y \in X$. 
Additionally, $X$ is \emph{Clifford} if the idempotents $xx^{-1}$ and $x^{-1}x$ are equal, for any $x\in X$, or, equivalently, the idempotents are central in the sense that they commute with every element in $X$. Furthermore, any Clifford semigroup is a \emph{union of groups.} More specifically, denoting by $\E(X)$ the set of all idempotents of $X$,
there exist a family $\{G_e\}_{e\in \E(X)}$ of disjoint groups and a family of group homomorphisms $\varphi_{f,e} : G_f\to G_e$, for all $e,f\in \E(X)$ with $e \leq f$, i.e., $ef=fe=e$, with $\varphi_{f,e}(y) = ey$, for any $y \in G_f$,
such that the product in $X$ can be expressed as $xy=  \varphi_{e, ef}\left(x\right)\varphi_{f, ef}\left(y\right)$, 
for all $x\in G_e$ and $y\in G_f$.

Let us observe that every Clifford semigroup $X$ gives rise to the following PE solutions 
\begin{align*}
\mathcal{I}(x,y)=(xy, y), \qquad \mathcal{F}(x,y)=\left(xy, yy^{-1}\right), \qquad \mathcal{E}(x,y)=(xy, e), 
\end{align*}
for all $x, y \in X$, where $e \in \E(X)$. These examples suggest studying the PE solution classes defined as follows.

\begin{definition}
A PE solution $s(x,y)=(xy, \theta_x(y))$  on a semigroup $X$ is said to be:
\begin{enumerate}
    \item  \emph{$\E\left(X\right)$-invariant} if $\theta_x(e) = \theta_x(f)$,
for all $x\in X$ and $e,f\in\E(X)$;
\item \emph{$\E(X)$-fixed} if $    \theta_x(e)=e$,
for all $x \in X$ and $e \in \E(X)$.
\end{enumerate}
\end{definition}
  In particular, $\mathcal{E}$ is $\E\left(X\right)$-invariant, while $\mathcal{I}$ and $\mathcal{F}$ are $\E\left(X\right)$-fixed on any Clifford semigroup $X$.

\begin{example}
Let $M =\{1,\,x,\,y\}$ be the commutative Clifford monoid defined by $x^2=x$, $y^2=x$, and $xy=y$, and $\gamma \in \End(M)$ the semigroup morphism given by $\gamma(1)=\gamma(x)=x$ and $\gamma(y)=y$. Then, the map $s(a,b)=(ab, \gamma(b))$ is an $\E(M)$-invariant PE solution on $M$.
\end{example}

Next, following \cite[Theorem 15]{MaPeSt23x}, we show how to construct all the $\E(X)$-invariant PE solutions on a Clifford semigroup $X$, by extending the description in \cref{teo_gruppi}. 
 First, recall that a congruence $\tau$ on $\E(X)$ is said to be \emph{normal} if  
\begin{align*}
    \forall \ e,f\in \E(X) \quad e \ \tau \ f \ \Longrightarrow \ \forall \ a\in X \quad a^{-1}ea\ \tau \ a^{-1}fa.
\end{align*}
Moreover, a subsemigroup $K$ of $X$ is \emph{normal} if $\E(X) \subseteq K$ and for all $a \in K$ it holds that $a^{-1}\in K$  and $a^{-1}Ka \subseteq K$. Thus, if $K$ is a normal subsemigroup of $X$, the pair $(K,\tau)$ is named a \emph{congruence pair} of $X$ if
\begin{align*}
  \forall\ a \in X,\ e\in \E(X)  \quad ae\in K \ \  \text{and} \ \ (e, a^{-1}a)\in \tau \ \Longrightarrow \ a\in K.
\end{align*}
Besides, recall that the restriction of a congruence $\rho$ in $X$ to $\E(X)$ is also a congruence on $\E(X)$,  called the \emph{trace} of $\rho$ and usually written as $\tr \rho$. Denote by $\Ker \rho$ the union of all the idempotent $\rho$-classes, namely,
\begin{align*}
    \Ker \rho = \bigcup_{e \in \E(X)} \rho e
\end{align*}
where $\rho e=\{x \in X \, \mid \, (x, e) \in \rho\}$, for all $e \in \E(X)$. The properties of $\rho$ can be described entirely in terms of $\Ker \rho$ and $\tr \rho$, as we show next (for more details, see \cite[Section 5.3]{Ho95}).

\begin{theorem}[cf. Theorem 5.3.3 in \cite{Ho95}]\label{theo_rho}
Let $X$ be an inverse semigroup. If $\rho$ is a congruence on $X$, then $(\Ker\rho, \tr\rho)$ is a congruence pair. Conversely, if $(K, \tau)$ is a congruence pair, then
\begin{align*}
    \rho_{(K, \tau)}=\lbrace (a,b) \in X \times X \, \mid \, \left(a^{-1}a, b^{-1}b\right) \in \tau, \, ab^{-1} \in K\rbrace
\end{align*}
is a congruence on $X$. Moreover, $\Ker \rho_{(K, \tau)}=K$, $\tr \rho_{(K, \tau)}=\tau$, and $\rho_{(\Ker \rho, \tr \rho)}=\rho$.
\end{theorem}

\begin{theorem}\label{propmu}
	Let $X$ be a Clifford semigroup, $\rho$ a congruence on $X$ such that $X/\rho$ is a group and $\mathcal{R}$ a system of representatives of $X/\rho$.  If $\mu: X \to \mathcal{R}$ is a map such that 
	\begin{equation*} \mu\left(xy\right)=\mu\left(x\right)\mu\left(x\right)^{-1}\mu\left(xy\right)
	\end{equation*}
	and $\left(\mu(x) , x\right) \in \rho$, for all $x,y \in X$,  then the map  $s(x,y)=\left(xy, \mu\left(x\right)^{-1}\mu\left(xy \right)\right)$,
	for all $x,y \in X$, is an $\E(X)$-invariant PE solution on $X$.
 
 Conversely, let $s(x,y)=(xy, \theta_x(y))$ be an $\E\left(X\right)$-invariant PE solution on $X$. Then, 
 there exists a congruence pair $\left(K, \tau\right)$, with 
 \[K=\{x\in X\, \mid \,\theta_e(x) \in \E(X), \,\text{for any} \ e \in \E(X)\} \text{ and } \tau=\E(X) \times \E(X),\] such that \begin{align*}
    \rho_{(K, \tau)}=\lbrace (x,y) \in X \times X \, \mid \, \left(x^{-1}x, y^{-1}y\right) \in \tau, \, xy^{-1} \in K\rbrace
\end{align*}
is a congruence on $X$ and
$\theta_e\left(X\right)$ is a system of representatives of the group $X/\rho_{\left(K, \tau\right)}$ and 
$\left(\theta_e\left(x\right), x \right) \in \rho_{\left(K, \tau\right)}$, for all $e \in \E(S)$ and $x\in X$. Moreover, the map $\theta_e$ satisfies 
\begin{equation*} \theta_e\left(xy\right)=\theta_e\left(x\right)\theta_e\left(x\right)^{-1}\theta_e\left(xy\right),
	\end{equation*}
 and $\theta_x(y)=\theta_e(x)^{-1}\theta_e(xy)$,
for all $x, y\in X$ and $e\in\E\left(X\right)$.
\end{theorem}

\begin{proposition}
 Let $X$ be a Clifford semigroup and 
 \[s(x,y)=(xy, \theta_x(y)), \: t(u,v)=(uv, \eta_u(v))\] two $\E(X)$-invariant PE solutions on $X$. Then, $s$ and $t$ are isomorphic if, and only if, there exists an isomorphism $\psi$ of $X$ such that $\psi\theta_e=\eta_e\psi$, for every $e \in \E(X)$.   
\end{proposition}

\noindent In other words,  $\psi$ sends the system of representatives $\theta_e(X)$ into the other one $\eta_e\left(\psi(X)\right)$.

With regards to the class of $\E(X)$-fixed PE solutions, a construction has been provided. Indeed, considering that every Clifford semigroup can be seen as a union of groups $G_e$, it is natural to contemplate whether it is possible to construct a global PE solution from PE solutions obtained in each of its groups.
\begin{proposition} \label{theo_costru}
Let $X$ be a Clifford semigroup and assume that, for every $e \in \E(X)$, $s^{[e]}(x,y)=\left(xy,\, \theta^{[e]}_x\left(y\right)\right)$ is a PE solution on $G_e$. Moreover,  for all $e, f\in \E(X)$, let $\epsilon_{e, f}: G_{e} \to G_{f}$ be maps such that $\epsilon_{e, f} = \varphi_{e, f}$ if $e\geq f$. Assume the following conditions are satisfied
\begin{align*}
   \theta^{[h]}_{\epsilon_{ef, h}(xy)} = \theta^{[h]}_{\epsilon_{e, h}(x)\epsilon_{f, h}(y)} \quad \text{and} \quad
    \epsilon_{f, h}\theta^{[f]}_{\epsilon_{e, f}(x)}(y) = \theta^{[h]}_{\epsilon_{e, h}(x)}\epsilon_{f, h}(y),
\end{align*}
for all $e, f, h\in \E(X)$ and $x\in G_{e}, y \in G_{f}$. Set $\theta_x(y) := \theta^{[f]}_{\epsilon_{e, f}(x)}(y),$ for all $x \in G_{e}, y \in G_{f}$. Then, the map $s(x,y)=(xy, \theta_x(y))$ is an $\E(X)$-fixed  PE solution on $X$.
\end{proposition}

\begin{example}
Let $X$ be a Clifford semigroup and $s^{[e]}\left(a,b\right) = \left(ab, \gamma^{[e]}\left(b\right)\right)$ a PE solution on the group $G_e$, with $\gamma^{[e]}$ an idempotent endomorphism of $G_e$, for every $e\in \E(X)$. Then, by choosing maps $\epsilon_{e,f}:G_e\to G_f$, for all $e,f\in \E(X)$ such that $\varphi_{e,f}\gamma^{[e]} = \gamma^{[f]}\varphi_{e,f}$ if $e\geq f$ and $\epsilon_{e,f}\left(x\right) := f$ otherwise,  one has that the map  $s(a,b)=\left(ab, \gamma^{[f]}(b)\right)$,
for all $a \in G_e$ and $b \in G_f$, is a PE solution on $X$.
\end{example}

Denoting by $K_e$ the kernel of each PE solution $s^{[e]}$ on $G_e$,  i.e., the normal subgroup $K_e=\{a \in G_e \, \mid \, \theta^{[e]}_{e}(a)=e\}$ of $G_e$ that we know to exist from \cref{teo_gruppi}, we have the following result. 

\begin{proposition}
    Let $X$ be a Clifford semigroup, $\epsilon_{e,f}(e)=f$, for all $e,f \in \E(X)$ with $e \leq f$, and $s$ be an $\E(X)$-fixed solution on $X$ constructed as in \cref{theo_costru}.   Then, $ K=\{a \in X \, \mid \, \forall\, e \in \E(X), \,  e \leq a, \ \theta_e(a)=aa^{-1}\} = \displaystyle \bigcup_{e\in \E(X)}\,K_{e}$.
\end{proposition}

However, in \cite{MaPeSt23x}, it is argued that this cannot be a description. In light of the previous discussion, the following question arises.
\begin{prob}
Find all the $\E(X)$-fixed PE solutions on a Clifford semigroup $X$.
\end{prob}

To conclude, we observe that not every PE solution on a Clifford $X$ lies in the class of $\E(X)$-invariant or  $\E(X)$-fixed PE solutions (see \cite[Example 4]{MaPeSt23x}). Thus, the next problem arises.
\begin{prob}
Find and study other classes of PE solutions on a Clifford semigroup $X$.  
\end{prob}

\section{Description of the involutive PE solutions}
In this section, we recall the description of the involutive PE solutions obtained by Colazzo, Jespers, and Kubat in \cite{CoJeKu20}.

 Initially, let $s(x, y)=(xy, \theta_x(y))$ be an involutive PE solution on a semigroup $X$. Then, beyond \eqref{p_one} and \eqref{p_two}, the following additional equalities
 \begin{align*}
     xy\theta_x(y)=x \quad \text{and} \quad \theta_{xy}\theta_x(y)=y
 \end{align*}
hold, for all $x, y \in X$. We give some immediate considerations.

\vspace{2mm}
\begin{remark}\hspace{1mm}
\begin{enumerate}
    \item If $G$ is a group, then the only bijective PE solution $s_G(x, y)=(xy, y)$ on $G$ is involutive if, and only if, $G$ is an elementary abelian $2-$group.
    \item If $L$ is a left zero semigroup and $s_L$ is an involutive PE solution on $L$, one has that the map $s=s_L\times s_G$ is an involutive PE solution on the left group $X = L \times G$. In addition, \cite[Theorem 3.2]{CoJeKu20} shows that the converse is true, namely, all the involutive PE solutions are of this type. It is mainly a consequence of \cite[Theorem I.1.27]{ClPr61}.
\end{enumerate}
\end{remark}
Therefore, from the considerations above, the problem of describing all the involutive PE solutions reduces to finding them on a left zero semigroup. It will turn out that all such solutions are
uniquely determined by two elementary abelian $2-$groups. Thus, in the following, we list some results and useful notions to give the main theorem.

First, following \cite[Section 4]{CoJeKu20}, given an involutive PE solution $s(x, y)=(xy, \theta_x(y))$ on a semigroup $X$, one can consider the following  congruence relation
\begin{align*}
   \forall x, y\in X \qquad  x \sim y \iff \theta_x=\theta_y.
\end{align*}
So,  setting $\overline{X}:=X/\sim$ and denoting by $\overline{x}$ the $\sim$-class of $x \in X$, one has that $\overline{X}$ is a left zero semigroup and the map $\Ret(s)\left(\overline{x}, \overline{y}\right)=\left(\overline{xy}, \overline{\theta}_{\overline{x}}\left(\overline{y}\right) \right)$ is an involutive PE solution on $\overline{X}$, called the \emph{retract} of $s$. Moreover, the PE solution $s$ is called \emph{irretractable} if $s=\Ret(s)$. One can show that $\Ret(s)$ is irretractable. 

The following is a complete description of irretractable
involutive PE solutions.

\begin{proposition}
    Let $(A, +)$ be an elementary abelian $2$-group. Then, $t_A(x, y) = (x, x + y)$ is an irretractable involutive PE solution on $A$. 
    
    Conversely, if  $s$ is an irretractable involutive PE solution on a semigroup $X$, then there exists
a structure $(X, +)$ of elementary abelian $2$-group on $X$ such that $s(x, y) =
(x, x + y)$, for all $x, y \in X$.
\end{proposition}

The following problem might be interesting to deal with.
\begin{prob}
   Study a possible retract relation for bijective (not necessarily involutive) PE solutions. 
\end{prob}

\begin{lemma} \cite[Proposition 5.1 - 5.3]{CoJeKu20}
    Let $(A, +)$ be an elementary abelian $2$-group , $t_A$ the irretractable involutive PE solution on $A$, $X$ a non-empty set, and $\sigma : A \to
\Sym(X)$ a map. Setting $S = X \times A$, we have that the map $\Ext^\sigma_X(t_A) : S \times S \to S \times S$ given by
$$\Ext^\sigma_X(t_A) ((x, a), (y, b)) = \left(\left(x, a\right), \left(\sigma_{a+b}\sigma^{-1}_{b}\left(y\right), a + b\right)\right),$$
for all $(x,a), (y,b) \in S$, is an involutive PE solution on $S$, called the \emph{extension of $t_A$ by $X$ and $\sigma$}. In addition, every involutive PE solution on a left zero semigroup can be so constructed.
\end{lemma}

\begin{theorem}\cite[Theorem 5.5 - 5.6]{CoJeKu20}
  Let $s(x,y)=(xy, \theta_x(y))$ be an involutive PE solution on a semigroup $S$. Then, there exist two elementary abelian $2$-groups $A$ and $G$, and a non-empty set $X$ such that $S$ may be identified with $X \times A \times G$ and 
  \begin{align*}
      s=\Ext^\sigma_X(t_A) \times s_G,
  \end{align*}
  for some $\sigma: G \to \Sym(X)$, where $t_A$ is the irretractable involutive PE solution defined on $A$ and $s_G$ is the bijective PE solution on $G$. Moreover, $\Ret(s)=t_A$.
  
  In addition, all the involutive PE solutions, up to isomorphism,  on a non-empty set $S$ are in bijective correspondence with decompositions of $S$ as a product of  $X \times A \times G$.
    \end{theorem}

Finally, \cite[Corollary 5.7]{CoJeKu20} determines how many involutive PE solutions there are. Namely, assume that
$S$ is a non-empty set of cardinality $2^n(2m+1)$, with $n,m \geq 0$. Then, there exists exactly $\binom{n+2}{2}$ involutive PE solutions on $S$.

\section{Idempotent PE solutions}
This section is devoted to recalling essential features of idempotent PE solutions. In particular, we mention the description of solutions on monoids having central idempotents, contained in \cite{Maz23}.

 First, a PE solution $s(x, y)=(xy, \theta_x(y))$ on a semigroup $X$ is idempotent if, and only if, the following equalities 
 \begin{align*}
     xy\theta_x(y)=xy \qquad \text{and} \qquad \theta_{xy}\theta_x(y)=\theta_x(y)
 \end{align*}
are satisfied, for all $x, y \in X$. 
Below we list properties, showing that idempotents of the semigroup $X$ play an important role. First, following \cite[p. 22]{ClPr61}, given a semigroup $X$ and $e \in \E(X)$, then $e$ is a \emph{left identity} (resp. \emph{right identity}) if $ex=x$ (resp. $xe=x$), for every $x \in X$, and the sets
\begin{align*}
    eX&=\{x \in X \, \mid \, ex=x\},\qquad Xe=\{x \in X \, \mid \, xe=x\}
\end{align*}
coincide respectively with the principal right and left ideals of $X$ generated by $e$. 

\enlargethispage{-2\baselineskip}

\begin{proposition}\label{a}
Let $X$ be a semigroup, $e \in \E(X)$, and $s(x,y)=(xy, \theta_x(y))$ an idempotent PE solution on $X$.
\begin{enumerate} 
    \item If $x \in Xe$, then
    \begin{enumerate}
    \item[a.] $x \in X\theta_x(e)$,
    \item[b.]  $\forall y \in X \quad \theta_y(x) \in X\theta_x(e) $,
    \item[c.] $\theta_e=\theta_e\theta_x$.
\end{enumerate}
\item If $x \in eX$, then
\begin{enumerate}
    \item[d.]  $\theta_e(x) \in \E(X)$, 
    \item[e.]   $x \in X \, \theta_e(x)$,
     \item[f.] $\forall y \in X \quad \theta_y(x) \in X \theta_e(x)$,
    \item[g.] $\theta_x$ is an idempotent map.
\end{enumerate}
\end{enumerate}
\end{proposition}

Clearly, \cref{a} can be applied if $X$ is a monoid. The following Corollary collects the consequences of this specialization.

\begin{corollary}\label{prop_sidemp_monoide}
If $M$ is a monoid and $s(x,y)=(xy, \theta_x(y))$ an idempotent PE solution on $M$, then the following hold for every $x \in M$:
\begin{enumerate}
    \item  $\theta_1(x) \in \E(X)$ and, in particular, $\theta_1(1)=1$,
\item $\theta_x=\theta_{\theta_1(x)}$, 
    \item  $\theta_1=\theta_1\theta_x$, 
    \item $\theta_x$ is idempotent.
\end{enumerate}
\end{corollary}

As a consequence of the definition, if $M$ is a cancellative monoid the unique idempotent PE solution on $X$ is the map $s(x,y)=(xy, 1)$; it belongs to the class of PE solutions discussed in Example \ref{exsKey}-$2$. On the other hand,  even considering non-cancellative monoids of small orders, one can note that among the PE solutions, several idempotent ones do not belong to the class of solutions in Example \ref{exsKey}-$2$. The following is an easy example in \cite[Appendix B]{tesi}.

\begin{example}
Let $M=\{1,a,b\}$ be the commutative monoid with identity $1$ and multiplication given by $a^2=a, ab=a,  b^2=1.$ Then, there are 3 idempotent solutions, up to isomorphism:
\begin{enumerate}
    \item $s(x, y)=(xy, 1)$;
    \item $r(x,y)=(xy, \gamma(y))$, with  $\gamma :M \to M$ defined by $\gamma(1)=\gamma(b)=1$ and $\gamma(a)=a$;
    \item $t(x,y)=(xy, \theta_x(y))$, with $\theta_x:M \to M$ the map given by $\theta_x(1)=1, \, \theta_x(a)=a$, for every $x \in M$,  and $\theta_1(b)=\theta_b(b)=1$ and $\theta_a(b)=b$.
\end{enumerate}
\end{example}

Note that, in general, the map $\theta_1$ is not a homomorphism.

\begin{proposition}\label{e_m_z_m}
Let $M$ be a monoid having central idempotents and $s(x,y)=(xy, \theta_x(y))$ an idempotent PE solution on $M$. Then, the map $\theta_1$ is an idempotent monoid homomorphism from $M$ to $\E(M)$. Moreover, denoting by
$\ker \theta_1= \{(x,y) \in M \times M \, \mid \, \theta_1(x)=\theta_1(y)\}$, 
the following hold
\begin{enumerate}
 \item $\theta_1(M)$ is a system of representatives of $M/\ker \theta_1$;
     \item $(\theta_x(y) ,y) \in \ker \theta_1 $, for all $x, y \in M$.
 \end{enumerate}
\end{proposition}

In light of the previous results, we describe the idempotent PE solutions on monoids having central idempotents (see \cite[Theorem 19]{Maz23}).

\begin{theorem}
 Let $M$ be a  monoid having central idempotents and $\mu$ an idempotent monoid homomorphism from $M$ to $\E(M)$ such that, for every $x \in M$, $\mu(x)=e_x$, with $e_x \in \E(M)$ a right identity for $x$.
 Moreover, let $\{ \theta_e: M \to M \, \mid \, e \in \im \mu\}$ be a family of maps such that $\theta_1=\mu$, and for all $e, f \in \im \mu$,
  \begin{align*}
     &\theta_e=\theta_e\theta_{ef},\\
     &\theta_e(xy)=\theta_e(x)\theta_h(y),\\
     &\theta_{ej}\theta_e(x)=\theta_e(x),
 \end{align*} 
 for all $x,y \in M$,  with $h=\mu(ex)$ and $j=\mu(x)$. Then,   $s(x,y)=(xy, \theta_{\mu(x)}(y))$ is an idempotent PE solution on $M$. Conversely, every idempotent PE solution on $M$ can be so constructed.
 \end{theorem}

\begin{prob}
    Describe idempotent PE solutions on other classes of semigroups.
\end{prob}

\section{Commutative and cocommutative PE solutions}
Other classes of PE solutions that can be studied are the commutative and the cocommutative ones. These classes were introduced in \cite{CaMaMi19} as an analogue of the commutative and the cocommutative multiplicative unitary operators, i.e., vector solutions defined on Hilbert spaces (for more details, see \cite[Definition 2.1]{BaSk93}).

\begin{definition}
	A PE solution $s: X \times X \to X \times X$ on a semigroup $X$ is said to be \emph{commutative} if $s_{12}s_{13}=s_{13}s_{12}$; \emph{cocommutative} if $s_{13}s_{23}=s_{23}s_{13}$.
\end{definition}
  It is a routine computation to check that a PE solution $s(x,y)=(xy, \theta_x(y))$ on a semigroup $X$ is commutative if and only if 
\begin{align*}
  xzy=xyz  \qquad \text{and}\qquad \theta_x=\theta_{xy}, 
\end{align*}
for all $x,y,z\in X$. 
On the other hand, $s$ is cocommutative if and only if 
\begin{align*}
    x\theta_y(z)=xz \qquad \text{and} \quad \theta_x\theta_{y}=\theta_y\theta_x,
\end{align*}
for all $x,y,z \in X$.
There exist PE solutions that are both commutative and cocommutative, such as the maps in $1.$ of Example \ref{exsKey}.  Moreover, according to \cite[Corollary 3.4]{CoJeKu20}, if $s$ is an involutive PE solution, then $s$ is both commutative and cocommutative. Besides, if $X$ is a monoid, the unique cocommutative PE solution is given by  $s(x,y)=(xy,y)$. All the commutative PE solutions on a monoid $X$ are described in \cite[Proposition 5]{Maz23}. The same results can be observed in the case of a commutative or cocommutative PE solutions on a Clifford semigroup $X$ (see \cite[Proposition 9]{MaPeSt23x}). We summarize below:

\begin{proposition}
Let $X$ be a monoid (or a Clifford semigroup). Then, a PE solution $s(x,y)=(xy, \theta_x(y))$ on $X$ is 
\begin{enumerate}
    \item commutative if and only if $X$ is a commutative monoid (or Clifford semigroup) and $\theta_x=\gamma$, with $\gamma \in \End(X)$, $\gamma^2=\gamma$, for every $x \in X$;
    \item cocommutative if and only if $\theta_x=\id_X$, for every $x \in X$.
\end{enumerate}

\end{proposition}

\begin{prob}
    Classify commutative or cocommutative PE solutions on other classes of semigroups.
\end{prob}

\subsubsection*{Acknowledgments}
  This work was partially supported by the Dipartimento di Matematica e Fisica ``Ennio De Giorgi'' - Università del Salento and by ``INdAM - GNSAGA Project" - CUP E53C22001930001. The author is a member of GNSAGA (INdAM) and the non-profit association ADV-AGTA. 

The author would like to thank Ivan Kaygorodov for inviting them to hold a seminar as part of the ``\emph{European Non-Associative Algebra Seminar}" and for giving them the possibility to share their knowledge on this topic. The author would also like to thank F. Catino, M.M. Miccoli, V. Pérez-Calabuig, and P. Stefanelli for their collaboration over the years.

The author thanks the referees for reading carefully their manuscript and for the suggestions that helped them to improve the paper.

{\small
    \def\cprime{$'$}
    
}

\EditInfo{January 1, 2024}{April 2, 2024}{David Towers and Ivan Kaygorodov}

\end{document}